\pdfoutput=1
\documentclass[11pt,a4paper]{article}
\usepackage[T1]{fontenc}
\usepackage[utf8]{inputenc}
\usepackage{lmodern}
\usepackage[margin=28mm]{geometry}
\usepackage{amsmath,amssymb,amsthm,mathtools}
\usepackage{microtype}
\usepackage[colorlinks=true,linkcolor=blue,urlcolor=blue,citecolor=blue]{hyperref}
\newtheorem{theorem}{Theorem}
\newtheorem{lemma}{Lemma}
\newtheorem{corollary}{Corollary}
\newcommand{\E}{\mathrm{E}}
\newcommand{\Pp}{\mathrm{P}}
\newcommand{\N}{\mathbb N}
\newcommand{\Z}{\mathbb Z}
\newcommand{\dist}{\operatorname{dist}}
\hypersetup{pdfauthor={Jan Vecer},pdftitle={A half-order remainder for the Chow-Robbins stopping boundary}}
\title{An $O(n^{-1/2})$ remainder for the\\ Chow--Robbins stopping boundary}
\author{Jan Vecer\thanks{Department of Probability and Mathematical Statistics,
Faculty of Mathematics and Physics, Charles University,
Sokolovsk\'a 83, 186 75 Prague 8, Czech Republic.
Email: \href{mailto:vecer@karlin.mff.cuni.cz}{\nolinkurl{vecer@karlin.mff.cuni.cz}}.
ORCID: \href{https://orcid.org/0000-0002-1501-4646}{0000-0002-1501-4646}.}}
\date{}
\begin{document}
\maketitle
\begin{abstract}
We sharpen the asymptotic approximation to the optimal stopping boundary
in the Chow--Robbins game, proving the inverse-square-root remainder
conjectured by Elton. The proof combines an exactly harmonic integral
for simple random walk with a Catalan first-passage renewal identity.
Averaging the stopping defect over one spatial period removes its leading
lattice dependence, and a coupling of states at equal depth transfers
this cancellation to the renewal sum. Consequently, the explicit integer
stopping rule can fail at no more than a constant times the square root
of the number of observation times.
\end{abstract}
\noindent\textbf{Keywords:} Chow--Robbins game; optimal stopping; random walk; stopping boundary; Catalan number; lattice correction.

\smallskip\noindent\textbf{2020 Mathematics Subject Classification:} Primary 60G40; secondary 60G50.

\section{Introduction and main result}
In the Chow--Robbins game, a player observes successive tosses of a fair
coin and chooses when to stop. Let $X_i=1$ for heads and $X_i=-1$ for
tails, with independent tosses, and let $S_n=\sum_{i=1}^nX_i$ be the
cumulative gain after $n$ tosses. Stopping at time $n\ge1$ gives the
average gain $S_n/n$. Equivalently, the player can be paid the fraction
of heads observed, which is $1/2+S_n/(2n)$ and gives the same optimal
rule. There is no fixed terminal horizon: the objective is to maximise
$\E(S_\tau/\tau)$ over almost surely finite stopping times $\tau\ge1$
based on the tosses seen so far. The question is how large the current
cumulative gain must be to make stopping preferable to continuing.

Chow and Robbins~\cite{ChowRobbins} established the existence and
structure of an optimal rule. Dvoretzky~\cite{Dvoretzky} developed
further structural properties, while Shepp~\cite{Shepp} and
Walker~\cite{Walker} identified the Brownian boundary and the leading
square-root asymptotics of its random-walk counterpart. The broader
connection between optimal stopping and free-boundary problems is
developed systematically by Peskir and Shiryaev~\cite{PeskirShiryaev}.

Approximating a discrete stopping boundary by a Brownian boundary
also requires understanding the error introduced by discrete
observations and spatial overshoots. Chernoff~\cite{Chernoff} and
Chernoff and Petkau~\cite{ChernoffPetkau} developed corrected
random-walk approximations in sequential stopping problems.
Siegmund~\cite{Siegmund} related diffusion corrections for boundary
crossing to ladder-height moments, and Lai, Yao, and
AitSahlia~\cite{LaiYaoAitSahlia} established a general continuity
correction for suitably initialized finite-horizon free-boundary
approximations. The present paper studies the finer remainder in
the asymptotic boundary of the infinite-horizon Chow--Robbins game.

Computational work on the game includes the experimental approach of
Medina and Zeilberger~\cite{MedinaZeilberger} and the rigorous computer
analysis of H\"aggstr\"om and W\"astlund~\cite{HaggstromWastlund}.
Christensen and Fischer~\cite{ChristensenFischer} obtained a Brownian
upper bound for the $S_n/n$ value function under a sub-Gaussian
increment assumption, together with a tight lower bound for the
Bernoulli walk. These bounds enable rigorous determination of
stopping and continuation states by backward induction. Their work
provides a complementary approach to the boundary expansion studied
here: it controls the discrete value function through comparison
with the continuous-time problem.

Elton~\cite[Theorem 1.1]{Elton} identified the term of order $n^{-1/4}$
in the real stopping boundary and proved an $O(n^{-7/24})$ remainder.
The present paper proves the $O(n^{-1/2})$ remainder conjectured by
Elton~\cite[(1.2)]{Elton}. We first state the result, then explain the
estimate that makes the sharper remainder possible.

\subsection{The boundary and the main result}
Suppose that the current cumulative gain is $u$ after $n\ge1$
observations. By independence of future tosses, the optimal expected
payoff from this state is
\[
 V(u,n)=\sup_{\tau}\E\frac{u+S_\tau}{n+\tau},
 \qquad u\in\mathbb R,\quad n\in\N,
\]
where here $S_0=0$ and $S_k$ denotes the simple symmetric random walk
of future increments. The supremum is over almost surely finite
nonnegative integer-valued stopping times; choosing $\tau=0$ means
stopping immediately and receiving $u/n$. Allowing real $u$ extends
the integer states of the coin-tossing game and lets us study the
boundary before rounding it to a reachable state.

We use the classical strictly increasing real boundary $\beta_n$:
stopping is optimal precisely when $u\ge\beta_n$, and its first crossing
is an optimal rule. These facts, including continuity in the real
spatial variable, are recalled in Elton~\cite[Section 1.1]{Elton}.
The corresponding Brownian boundary has the form $\alpha\sqrt n$,
where $\alpha\in(0,1)$ is determined by
\begin{equation}\label{alpha}
 \alpha=(1-\alpha^2)\int_0^\infty \mathrm{e}^{\alpha s-s^2/2}\,\mathrm{d}s.
\end{equation}
Numerically, $\alpha=0.839923675692\ldots$. The only quantitative
boundary estimate imported into the proof is
\begin{equation}\label{coarse}
 \beta_n=\alpha\sqrt n-\frac12+O(n^{-1/4}),
\end{equation}
which follows from Elton~\cite[Lemma 4.4(a)]{Elton}. In particular,
the argument below does not assume the coefficient or the $7/24$
remainder from his main theorem.

\begin{theorem}\label{main}
As $n\to\infty$ through the positive integers,
\begin{equation}\label{result}
 \boxed{\displaystyle
 \beta_n=\alpha\sqrt n-\frac12
 -2\zeta(-1/2)\sqrt{\frac\alpha\pi}\,n^{-1/4}
 +O(n^{-1/2}).}
\end{equation}
\end{theorem}
The coefficient in~\eqref{result} is exactly the one obtained by
Elton~\cite[Theorem 1.1]{Elton}; the new conclusion is the remainder
conjectured in Elton~\cite[(1.2)]{Elton}.

\noindent\emph{Organisation of the proof.}
The proof of Theorem~\ref{main} occupies Sections 2--7.
Sections 2--6 establish the required estimates and reduce the boundary
equation to the scalar identity~\eqref{scalar}; Section 7 evaluates
this identity and completes the proof. The main steps are outlined
in Section~\ref{proof-outline} below.

Throughout, $C$ denotes a finite constant that may change
from line to line. All estimates asserted to be uniform use constants
independent of the indicated state or summation index.

\subsection{Where the sharper remainder comes from}\label{proof-outline}
Elton's proof uses a truncated backward-induction tree with Catalan
weights, Brownian value bounds, and feedback estimates for the values
at its leaves. At the final stage, the truncation parameter $J$ gives
two error terms,
\[
 O(J^{-1/2}n^{-1/4})+O(J^{5/2}n^{-1/2}).
\]
Balancing them with $J\asymp n^{1/12}$ gives the exponent $7/24$;
see Elton~\cite[Section 6]{Elton}. Our proof retains the Catalan
first-passage structure and controls the aggregate error through
spatial averaging. The argument has three steps.

\begin{enumerate}
\renewcommand{\labelenumi}{(\roman{enumi})}
\item \emph{An exact harmonic representation.}
Sections 2--3 construct a lattice analogue $H$ of Shepp's Brownian
continuation-value integral. Replacing
$\mathrm{e}^{-\lambda^2n/2}$ by $(\cosh\lambda)^{-n}$ makes the
integral exactly harmonic for the random walk. Optional stopping
expresses the defect $H-V$ through its terminal value. After
normalisation, its leading dependence on the distance below the
Brownian boundary has a periodic part with period one.

\item \emph{Cancellation of the lattice phase.}
Sections 4--5 average the defect over one spatial unit. This average
preserves harmonicity and replaces the leading periodic part by its
mean, $1/12$ (Lemma~\ref{meanlemma}). A coupling of states at equal
depth below the Brownian boundary transfers the fixed-time spatial
average to the renewal curve (Lemma~\ref{alignment}). Summing by
unit-depth cells then gives an $O(n^{-1/2})$ bound for the weighted
stopping-phase error (Lemma~\ref{averagederror}). This is the decisive
estimate: it preserves cancellation in the signed aggregate that
separate bounds on the individual leaf errors would discard.

\item \emph{Reduction to a scalar renewal sum.}
The exact Catalan renewal identity in Section 5 converts these defect
estimates into an equation for the boundary. Section 6 controls the
signed error from replacing the square-root phase along the renewal
curve by its linearisation. Section 7 evaluates the remaining scalar sum, recovering the
coefficient involving $\zeta(-1/2)$ with an $O(n^{-1/2})$ remainder.
\end{enumerate}
Thus the improvement rests on the averaged phase estimate and its
transfer to the renewal sum. It does not require an approximation to
the phase-dependent moving-boundary hitting law.

\subsection{The integer stopping rule}
The real-boundary estimate yields an explicit rule for integer states
outside a set of times of small density.

\begin{corollary}\label{integer-rule}
Put
\[
 c=-2\zeta(-1/2)\sqrt{\alpha/\pi},\qquad
 a(x)=\alpha\sqrt{x}-\frac12+cx^{-1/4},\qquad
 k_n=\lceil\beta_n\rceil.
\]
There is a set $E\subseteq\N$ satisfying
\[
 \#(E\cap\{1,\ldots,N\})=O(N^{1/2})
\]
such that $k_n=\lceil a(n)\rceil$ whenever $n\notin E$.
In particular, the proportion of exceptional times up to $N$ is
$O(N^{-1/2})$.
\end{corollary}
\begin{proof}
By Theorem~\ref{main}, there are $K>0$ and $n_0$ such that
$|\beta_n-a(n)|\le Kn^{-1/2}$ for $n\ge n_0$. Include the finitely
many smaller times in $E$, and set
\[
 E\cap[n_0,\infty)
 =\{n\ge n_0:\dist(a(n),\Z)\le Kn^{-1/2}\}.
\]
Outside this set the closed interval between $a(n)$ and $\beta_n$
contains no integer, so their ceilings agree.

For large $M$, throughout $[M,2M]$ we have
$a'(x)\asymp M^{-1/2}$. The range of $a$ on this interval has length
$O(M^{1/2})$. Thus only $O(M^{1/2})$ integers $j$ can satisfy
$|a(x)-j|\le KM^{-1/2}$ for some $x\in[M,2M]$.
For each such $j$, the set of these $x$ is an interval of length
$O(1)$, by the lower bound on $a'$. It contains
$O(1)$ integers. Consequently
$\#(E\cap[M,2M])=O(M^{1/2})$; summing over dyadic intervals proves
the assertion.
\end{proof}

\noindent\emph{Parity convention.}
The threshold $k_n$ is the canonical threshold on all integer spatial
states. Since $S_n\equiv n\pmod 2$, the least reachable stopping
height is
\[
 k_n^{\mathrm{reach}}
 =n+2\left\lceil\frac{\beta_n-n}{2}\right\rceil.
\]
For $n\notin E$ it is equally given by
$n+2\lceil(a(n)-n)/2\rceil$.

\section{An exact harmonic integral}
Write $A=1-\alpha^2$ and define, for $u<n$,
\begin{equation}\label{H}
 H(u,n)=A\int_0^\infty \mathrm{e}^{\lambda u}(\cosh\lambda)^{-n}\,\mathrm{d}\lambda.
\end{equation}
The domain $u<n$ is preserved by a random-walk step. Directly under the
integral,
\begin{equation}\label{harmonic}
 \frac12\{H(u-1,n+1)+H(u+1,n+1)\}=H(u,n).
\end{equation}
For each fixed $K$, uniformly for $u\le\alpha\sqrt n+K$ and large $n$,
\begin{equation}\label{Hbound}
 0\le H(u,n)\le C_K n^{-1/2}.
\end{equation}
Indeed, on $0\le\lambda\le1$ use $\log\cosh\lambda\ge\lambda^2/4$
and rescale $s=\lambda\sqrt n$. On $\lambda\ge1$, use
$\log\cosh\lambda\ge\lambda-\log2$. The integral over this latter
range is at most
$\mathrm{e}^{n\log2-(n-u)}/(n-u)$, exponentially small uniformly in the stated
region.

\begin{lemma}\label{Hex}
Let $\kappa=(4+\alpha^2)/6$. For $d$ in any fixed compact interval,
\begin{equation}\label{Hexeq}
 H(\alpha\sqrt n-d,n)-\frac{\alpha\sqrt n-d}{n}
 =\alpha n^{-3/2}(d^2+\kappa)+O(n^{-2}),
\end{equation}
uniformly in $d$.
\end{lemma}
\begin{proof}
Put $I(z)=\int_0^\infty \mathrm{e}^{zs-s^2/2}\,\mathrm{d}s$. Integration by parts gives
$I'(z)=zI(z)+1$. Together with~\eqref{alpha}, this yields
\[
 AI(\alpha)=\alpha,\quad AI'(\alpha)=1,\quad
 AI''(\alpha)=2\alpha,\quad
 AI^{(4)}(\alpha)=2\alpha(4+\alpha^2).
\]
After $s=\lambda\sqrt n$, expand the integrand using
\[
 -n\log\cosh(s/\sqrt n)
 =-s^2/2+s^4/(12n)+O(s^6/n^2).
\]
For example, on $s\le n^{1/12}$ the exponential expansion has an
integrable error bounded by a constant times
$n^{-2}(s^6+s^8)\mathrm{e}^{\alpha s-cs^2}$, with $c>0$. Outside that window
the same quadratic lower bound as in~\eqref{Hbound} gives an
exponentially small integral, and the $\lambda\ge1$ tail was already
controlled. Expanding $\mathrm{e}^{-ds/\sqrt n}$ uniformly for bounded $d$ gives
\[
 H(\alpha\sqrt n-d,n)
 =\alpha n^{-1/2}-d n^{-1}
 +\alpha\left(d^2+\frac{4+\alpha^2}{6}\right)n^{-3/2}
 +O(n^{-2}).
\]
Subtracting the immediate payoff proves the lemma.
\end{proof}

\section{The stopping defect and its lattice phase}
Define $D(u,n)=H(u,n)-V(u,n)$. For a continuation state, let $\tau$ be
the optimal first crossing of the boundary. This time is almost surely
finite, either by the classical stopping theorem or by the law of the
iterated logarithm and~\eqref{coarse}. Before crossing, the position is
below $\beta_{n+k}$; at crossing, its overshoot over the previous
boundary is at most one. Consequently~\eqref{Hbound} bounds the stopped
$H$-martingale uniformly, and dominated convergence gives
\begin{equation}\label{Drep}
 D(u,n)=\E\left[H(u+S_\tau,n+\tau)
       -\frac{u+S_\tau}{n+\tau}\right].
\end{equation}
Put $t=n+\tau$ and $d_\tau=\alpha\sqrt t-u-S_\tau$.
At the crossing, $u+S_\tau\ge\beta_t$ while
$u+S_\tau<\beta_{t-1}+1$. Thus~\eqref{coarse} implies
\begin{equation}\label{overshoot}
 -\frac12-Ct^{-1/4}\le d_\tau\le\frac12+Ct^{-1/4}.
\end{equation}
Lemma~\ref{Hex} and~\eqref{Drep} now prove the uniform bound
\begin{equation}\label{Dbound}
 |D(u,n)|\le Cn^{-3/2}
 \quad\text{for every continuation state and all large }n.
\end{equation}
The bound remains valid in the bounded strip of stopping states used
below, by Lemma~\ref{Hex} directly.

Let $h(d)=\dist(d,\Z)^2$. This function is bounded, 1-periodic and
globally 1-Lipschitz. Write $d=\alpha\sqrt n-u$. Since $S_\tau$ is an
integer,
\[
 h(d_\tau)=h\bigl(d+\alpha(\sqrt{n+\tau}-\sqrt n)\bigr).
\]
Equation~\eqref{overshoot} gives
$|d_\tau^2-h(d_\tau)|\le C(n+\tau)^{-1/4}$. The Lipschitz bound,
Lemma~\ref{Hex}, and
$|1-(n/(n+\tau))^{3/2}|\le C\min(\tau/n,1)$ consequently imply
\begin{equation}\label{phasepre}
 \left|\frac{n^{3/2}}\alpha D(u,n)-\kappa-h(d)\right|
 \le C\left\{n^{-1/4}+\E\min(\tau/\sqrt n,1)\right\}.
\end{equation}

For integer $0\le s\le\sqrt n$, the boundary rise is uniformly bounded, so
$\beta_{n+s}-u\le d+C$. Choose a positive integer $\ell$ above all
these levels, with $\ell\le d_++C$, where $d_+=\max(d,0)$.
If the walk reaches $\ell$ by time $s$, it has already crossed the
optimal boundary. The reflection principle and the usual binomial
mass bound give, for the first passage $T_\ell$ to $\ell$,
\[
 \Pp(\tau>s)\le\Pp(T_\ell>s)
 \le\min\left(1,\frac{C\ell}{\sqrt{s+1}}\right),
 \qquad 0\le s\le\sqrt n.
\]
For clarity, the last bound also follows by reflecting paths that hit
$\ell$: the surviving probability is bounded by the probability that
$S_{\lfloor s\rfloor}$ lies in an interval of $O(\ell)$ lattice sites,
each having probability at most $C/\sqrt{s+1}$.
Extending the tail bound to real $s$ by integer rounding, integrating it,
and using~\eqref{Dbound} for the trivial cap gives
\begin{lemma}\label{phase}
Uniformly over continuation states with large $n$,
\begin{equation}\label{phaseeq}
 \frac{n^{3/2}}\alpha D(u,n)
 =\kappa+h(d)+O\!\left(\min\{1,(1+d_+)n^{-1/4}\}\right).
\end{equation}
\end{lemma}
In a stopping state the corresponding leading term is
$\kappa+d^2$. Wrapping $d$ modulo one there is generally incorrect.

\section{Spatial averages and comparison of equal-depth states}
Write
\[
 F(u,t)=\frac{t^{3/2}}\alpha D(u,t),\qquad
 e_t(d)=F(\alpha\sqrt t-d,t)-\kappa-h(d).
\]
The following estimates use only the preliminary localization
\eqref{coarse} and Lemmas~\ref{Hex} and~\ref{phase}.

\begin{lemma}[Spatial regularity]\label{regularity}
For all sufficiently large $t$, the function $F(\cdot,t)$ is uniformly
bounded on the continuation region and is $C$-Lipschitz on
$(-\infty,\beta_t-1]$, with $C$ independent of $t$.
\end{lemma}
\begin{proof}
Boundedness follows from~\eqref{Dbound}. Differentiating the harmonic
integral twice gives, uniformly for $u\le\beta_t+1$,
\[
 0\le H_{uu}(u,t)\le Ct^{-3/2}.
\]
Indeed the integrand is increasing in $u$, and the rescaling used
in~\eqref{Hbound}, with the extra factor $\lambda^2$, gives this bound
at the upper endpoint. The value $V(\cdot,t)$ is convex, being a
supremum of affine functions. Its one-sided derivatives at $u$ lie
between $V(u,t)-V(u-1,t)$ and $V(u+1,t)-V(u,t)$. For
$u\le\beta_t-1$, all three defects are $O(t^{-3/2})$.
Comparison with the corresponding secants of $H$ therefore gives
\[
 |V'_\pm(u,t)-H_u(u,t)|\le Ct^{-3/2}.
\]
Local absolute continuity of the convex function $V$ now proves the
claimed Lipschitz bound.
\end{proof}

The argument uses one-sided derivatives and convexity, and does not
assume differentiability of the discrete value function.
This distinction is relevant to the present problem: Christensen and
Fischer~\cite{ChristensenFischerRegularity} show that discrete-time
stopping values can fail to be spatially differentiable on a dense
subset of the continuation region, with the Chow--Robbins game as
an example.

\begin{lemma}[Discount loss and spatial phase mean]\label{meanlemma}
For sufficiently large $t$, let $u<\beta_t$, put
$d=\alpha\sqrt t-u$, and let $\tau$ be its
optimal first-crossing time. Then
\begin{equation}\label{discountloss}
 \E\left[1-\sqrt{\frac t{t+\tau}}\right]
 \le C\min\left\{1,\frac{d+1}{\sqrt t}\right\}.
\end{equation}
Moreover,
\begin{equation}\label{phasemean}
 \left|\int_d^{d+1}e_t(y)\,\mathrm{d}y\right|
 \le C\left[t^{-1/4}+
       \min\left\{1,\frac{d+1}{\sqrt t}\right\}\right].
\end{equation}
\end{lemma}
\begin{proof}
Put $T=t+\tau$ and $X=u+S_\tau$. The crossing bound implies
$X\le\alpha\sqrt T+C$. Optimality and the immediate stopping payoff
give
\[
 \frac{\alpha\sqrt t-d}{t}=\frac ut
 \le V(u,t)=\E\frac XT
 \le\alpha\E T^{-1/2}+\frac Ct.
\]
Rearrangement proves~\eqref{discountloss}, including its trivial cap.
Since $1-z^3\le3(1-z)$ for $0\le z\le1$, it also gives
\begin{equation}\label{discountthree}
 1-\E(t/T)^{3/2}
 \le C\min\left\{1,\frac{d+1}{\sqrt t}\right\}.
\end{equation}

Consider the spatial average
\[
 \overline D(u,t)=\int_0^1D(u-r,t)\,\mathrm{d}r.
\]
It is harmonic whenever $u<\beta_t$, since the whole interval
$[u-1,u]$ lies in continuation. Stop its martingale when the upper
endpoint first crosses. The uniform defect bound, including the
bounded stopping strip at the terminal time, bounds the stopped
average by $Ct^{-3/2}$. Optional stopping therefore gives
\[
 \overline D(u,t)=\E\overline D(X,T).
\]
The terminal depth $d_\tau=\alpha\sqrt T-X$ satisfies
\eqref{overshoot}. For $0\le r\le1$, Lemma~\ref{phase} applies
to $X-r$ if it is continuing, and Lemma~\ref{Hex} applies if it
is stopping. A stopping depth $d_\tau+r$ lies in
$[-1/2-CT^{-1/4},1/2+CT^{-1/4}]$; hence its square differs from
$h(d_\tau+r)$ by $O(T^{-1/4})$. In both cases,
\[
 F(X-r,T)=\kappa+h(d_\tau+r)+O(T^{-1/4}),
\]
uniformly in $r$. The mean of $h$ over every unit interval is
$1/12$, so
\[
 \frac{t^{3/2}}\alpha\overline D(u,t)
 =\left(\kappa+\frac1{12}\right)\E(t/T)^{3/2}
  +O(t^{-1/4}).
\]
Subtracting $\kappa+1/12$ and using~\eqref{discountthree} proves
\eqref{phasemean}.
\end{proof}

\begin{lemma}[Comparison of equal-depth states]\label{alignment}
Fix $C_0<\infty$. Let $t,s$ be integers with
$0\le s\le C_0\sqrt t$, put
$r=\alpha(\sqrt{t+s}-\sqrt t)$, and suppose
$d=\alpha\sqrt t-x\ge2$. There is a finite constant $C_{C_0}$,
depending on $C_0$, such that uniformly over these states,
\begin{equation}\label{alignmentbound}
 |F(x+r,t+s)-F(x,t)|
 \le C_{C_0}\left[t^{-1/4}+
       \min\left\{1,\frac{d+1}{\sqrt t}\right\}\right].
\end{equation}
\end{lemma}
\begin{proof}
All constants in this proof may depend on $C_0$.
The two initial states are continuing for large $t$ and have the same
depth below their Brownian boundaries. Couple them using the same
random walk. Let their optimal crossing times be $\tau_A,\tau_B$,
and set $\sigma=\min(\tau_A,\tau_B)$ and $T=t+\sigma$.
Both defect processes are harmonic before $\sigma$ and uniformly
bounded when stopped there. At $\sigma$ at least one state has just
crossed. Its depth is in the strip~\eqref{overshoot}. The other
state either also crossed or is continuing, with depth at least
$1/2-CT^{-1/4}$. The two depths differ by at most $r=O(1)$.
Thus both terminal depths lie in a fixed compact interval with lower
endpoint $-1/2-CT^{-1/4}$.

Denote those depths by $d_A,d_B$. The stopping expansion and
Lemma~\ref{phase}, exactly as in the preceding proof, give
\[
 F(x+S_\sigma,T)=\kappa+h(d_A)+O(t^{-1/4}),\qquad
 F(x+r+S_\sigma,T+s)=\kappa+h(d_B)+O(t^{-1/4}).
\]
Their phase difference has the exact expression
\[
 d_B-d_A
 =r\left[\frac{\sqrt{t+s}+\sqrt t}{\sqrt{T+s}+\sqrt T}-1\right].
\]
The fraction lies between $\sqrt{t/T}$ and $1$: the lower bound
follows by comparing each denominator term with $\sqrt{T/t}$ times
its numerator term. Hence
\[
 |d_B-d_A|\le r(1-\sqrt{t/T}).
\]
Since $\sigma\le\tau_A$, Lemma~\ref{meanlemma} bounds the
expectation of the right side by
$C\min\{1,(d+1)/\sqrt t\}$. Also, uniformly in $\sigma$,
\[
 \left|\left(\frac{t+s}{T+s}\right)^{3/2}
              -\left(\frac tT\right)^{3/2}\right|
 \le C\frac st=O(t^{-1/2}).
\]
Use optional stopping for each defect process, insert the terminal
expansions, and use boundedness and Lipschitz continuity of $h$.
These estimates prove~\eqref{alignmentbound}.
\end{proof}

\section{Renewal at the real boundary}
Fix $u=\beta_n$ and write $\delta=\alpha\sqrt n-u$; thus
$\delta=1/2+O(n^{-1/4})$. Let $T_1$ be first passage to $+1$.
Its Catalan distribution is
\begin{equation}\label{Catalan}
 \Pp(T_1=2m+1)=p_m=2^{-2m-1}\frac{1}{m+1}\binom{2m}{m},
 \quad m\ge0.
\end{equation}
In particular $\sum p_m=1$, $p_m\le C(m+1)^{-3/2}$, and
\begin{equation}\label{pm}
 p_m=\frac{1}{2\sqrt\pi}m^{-3/2}+O(m^{-5/2}),\qquad m\ge1.
\end{equation}

The value is $1/n$-Lipschitz in $u$, since each stopped payoff has this
Lipschitz bound. Letting a continuation state tend to $\beta_n$ from
below therefore proves equality with its one-step continuation value
at the boundary. Before $T_1$, all subsequent positions are at most
$u$, and the strict increase of $\beta$ puts them in the continuation
region. Iterating Bellman equality yields
\[
 V(u,n)=\E V(u+1,n+T_1).
\]
The limiting step is justified, for example, by $0\le V\le1$ on the
reachable domain $u<n$. For each deterministic $N$, we have
$V(u,n)\ge\E[(u+S_N)/(n+N)]=u/(n+N)$; letting $N\to\infty$
proves the lower bound. The upper bound follows pathwise from
$u+S_k<n+k$. The analogous $H$ identity follows from its bounded
stopped martingale. Subtracting gives the exact renewal identity
\begin{equation}\label{renewal}
 D(u,n)=\sum_{m\ge0}p_mD(u+1,n+2m+1).
\end{equation}

Use the defect profile
\begin{equation}\label{profile}
 f(d)=\begin{cases}d^2,&d\le1/2,\\h(d),&d\ge1/2.\end{cases}
\end{equation}
Only its restriction to $d\ge-0.6$ is needed; there it is bounded and
globally Lipschitz. At a leaf let
$t_m=n+2m+1$ and $d_m=\alpha\sqrt{t_m}-u-1$.
For $m\le c_0\sqrt n$, with $c_0=1/(4\alpha)$, we have
$d_m\le-1/4+o(1)$. These leaves are stopping states and their
normalized expansion from Lemma~\ref{Hex} has error $O(n^{-1/2})$.
They include the early leaves with $d_m<-1/2$, for which replacing
$f$ by $h$ would change the leading coefficient.

For later leaves apply Lemma~\ref{phase} in continuation and
Lemma~\ref{Hex} in stopping. If a stopping leaf falls in
$1/2<d_m\le1/2+Ct_m^{-1/4}$, the difference $d_m^2-h(d_m)$ is
$O(t_m^{-1/4})$ and is absorbed in the same error. A continuing leaf
with $d_m<1/2$ is in a bounded neighborhood of $1/2$, where $f=h$.

\begin{lemma}[The averaged stopping-phase error]\label{averagederror}
With $u=\beta_n$, $t_m=n+2m+1$, and
$d_m=\alpha\sqrt{t_m}-u-1$, we have
\begin{equation}\label{averagedbound}
 \sum_{m\ge0}p_m\left(\frac n{t_m}\right)^{3/2}
 [F(u+1,t_m)-\kappa-f(d_m)]=O(n^{-1/2}).
\end{equation}
\end{lemma}
\begin{proof}
The early stopping leaves $m\le c_0\sqrt n$ contribute
$O(n^{-1/2})$ by Lemma~\ref{Hex}. For the remaining leaves with
$d_m<2$, the preceding classification gives pointwise error
$O(n^{-1/4})$; their total Catalan mass is $O(n^{-1/4})$.
Leaves with $m>n$ contribute $O(n^{-1/2})$ by boundedness of
the normalized defect and the profile. It remains to treat
$m\le n$ and $d_m\ge2$, where all leaves are continuing and $f=h$.

Put $x=u+1$ and $w_m=p_m(n/t_m)^{3/2}$. Partition the remaining
indices into unit-depth cells
\[
 I_j=\{m\le n:j\le d_m<j+1\},\qquad j\ge2.
\]
The last cell may be incomplete; it contains $O(\sqrt n)$ indices
with $m\asymp n$ and therefore has mass $O(n^{-1})$, which is
harmless. All other cells span one full phase period up to the mesh
error at their endpoints.

Write $q=\alpha\sqrt n-u-1=-1/2+O(n^{-1/4})$. The identity
\[
 2m+1=\frac{2\sqrt n}{\alpha}(d_m-q)
                 +\frac{(d_m-q)^2}{\alpha^2}
\]
shows that, on a cell with $m\le n$,
\begin{equation}\label{cellscales}
 j\le C\sqrt n,\quad m\asymp j\sqrt n,\quad
 \Delta_m:=d_{m+1}-d_m\asymp n^{-1/2},\quad
 W_j:=\sum_{m\in I_j}w_m\le Cn^{-1/4}j^{-3/2}.
\end{equation}
Let $m_j$ be the first index in a full cell and set $v_j=t_{m_j}$.
Define, for $j\le y\le j+1$,
\[
 g_j(y)=F(\alpha\sqrt{v_j}-y,v_j)-\kappa-h(y).
\]
Lemma~\ref{alignment}, applied to equal-depth states at times
$v_j$ and $t_m$, gives
\begin{equation}\label{celltransfer}
 F(x,t_m)-\kappa-h(d_m)
 =g_j(d_m)+O\left(n^{-1/4}+\frac{j+1}{\sqrt n}\right).
\end{equation}
Indeed $t_m-v_j=O(\sqrt n)$ throughout the cell. The fixed-time
phase mean, the spatial regularity lemma, and Lemma~\ref{phase}
give respectively
\begin{align}
 \left|\int_j^{j+1}g_j(y)\,\mathrm{d}y\right|
   &\le C\left(n^{-1/4}+\frac{j+1}{\sqrt n}\right),
       \label{cellmean}\\
 \operatorname{Lip}(g_j)&\le C,\qquad
 \sup_{[j,j+1]}|g_j|\le C A_j,
 \quad A_j=\min\{1,jn^{-1/4}\}.\label{cellamplitude}
\end{align}
The spatial arguments used here lie below $\beta_{v_j}-1$ for
large $n$, since $j\ge2$.

For completeness, the variation of the quadrature density is also
controlled. Put $a_m=w_m/\Delta_m$. The exact relations
\[
 \frac{p_{m+1}}{p_m}=\frac{2m+1}{2m+4},\qquad
 \Delta_m=\frac{2\alpha}{\sqrt{t_m+2}+\sqrt{t_m}}
\]
imply $|\log(a_{m+1}/a_m)|\le C((m+1)^{-1}+n^{-1})$.
There are $O(\sqrt n)$ indices in a cell, so
\[
 \sup_{m\in I_j}\left|\frac{a_m}{a_{m_j}}-1\right|
 \le C(j^{-1}+n^{-1/2}).
\]
Since $\sum_{m\in I_j}\Delta_m=1+O(n^{-1/2})$,
we have $a_{m_j}\asymp W_j$. The Lipschitz bound gives the
unit-cell quadrature estimate
\[
 \sum_{m\in I_j}\Delta_m g_j(d_m)
 =\int_j^{j+1}g_j(y)\,\mathrm{d}y+O(n^{-1/2});
\]
Here the last integration interval is truncated at $j+1$; its excess
and the initial endpoint gap are each charged to $O(n^{-1/2})$.
It follows from~\eqref{cellmean}--\eqref{cellamplitude} that
\begin{align*}
 \left|\sum_{m\in I_j}w_m g_j(d_m)\right|
 \le C W_j\left[
 n^{-1/4}+\frac{j+1}{\sqrt n}+n^{-1/2}
 +(j^{-1}+n^{-1/2})A_j\right].
\end{align*}
The transfer error in~\eqref{celltransfer} satisfies the same bound.
Using~\eqref{cellscales}, the mean and transfer terms sum to
\[
 Cn^{-1/4}\sum_{2\le j\le C\sqrt n}j^{-3/2}
       \left(n^{-1/4}+\frac{j+1}{\sqrt n}\right)
 =O(n^{-1/2}).
\]
The weight-variation terms sum to
\[
 Cn^{-1/4}\sum_{j\ge2}j^{-3/2}
       (j^{-1}+n^{-1/2})\min\{1,jn^{-1/4}\}
 =O(n^{-1/2}),
\]
as follows by splitting at $j=n^{1/4}$; the mesh terms are smaller.
Combining the cells, the initial leaves, and the tail proves
\eqref{averagedbound}.
\end{proof}

At the root, Lemma~\ref{Hex} gives
$F(u,n)=\delta^2+\kappa+O(n^{-1/2})$. Combining the exact renewal
identity with Lemma~\ref{averagederror} yields
\begin{equation}\label{renewalapprox}
 \delta^2+\kappa=
 \sum_{m\ge0}p_m(n/t_m)^{3/2}\{f(d_m)+\kappa\}
 +O(n^{-1/2}).
\end{equation}
Replacing $(n/t_m)^{3/2}$ by $1$ costs $O(n^{-1/2})$, because
\begin{equation}\label{factorremoval}
 \sum_{m\ge0}p_m\min((m+1)/n,1)=O(n^{-1/2}).
\end{equation}

\section{Signed curvature replacement}
\begin{lemma}[Signed curvature replacement]\label{curvaturelemma}
Let $q$ range over a compact interval and let $f$ be bounded and
Lipschitz on a half-line containing all $q+y$, $y\ge0$. Then, uniformly
in $q$,
\begin{align}\label{curvaturesigned}
\sum_{m\ge0}p_m\biggl[
f\!\left(q+\alpha(\sqrt{n+2m+1}-\sqrt n)\right)
-f\!\left(q+\frac{\alpha m}{\sqrt n}\right)
\biggr]=O(n^{-1/2}).
\end{align}
\end{lemma}
\begin{proof}
Write $B=\sqrt n$, $\eta=\alpha/B$, and
$g_B(t)=\alpha(\sqrt{B^2+2t}-B)$. Removing the $+1$ inside the
square root costs $O(B^{-1})$, since the change of argument is at
most $\alpha/(2B)$ and $\sum p_m=1$.

For either $a_B(t)=g_B(t)$ or $a_B(t)=\eta t$, set
\[
 z_B(t)=t^{-3/2}[f(q+a_B(t))-f(q)].
\]
Both choices satisfy $a_B(t)\le\alpha t/B$ and
$0\le a_B'(t)\le\alpha/B$. Therefore
\[
 |f(q+a_B(t))-f(q)|\le C\min(t/B,1).
\]
The Catalan asymptotic gives
\[
 \sum_{m\ge1}p_m[f(q+a_B(m))-f(q)]
 =\frac1{2\sqrt\pi}\sum_{m\ge1}z_B(m)+O(B^{-1}),
\]
because $\sum_{m\ge1}m^{-5/2}\min(m/B,1)=O(B^{-1})$.
Moreover, almost everywhere,
\[
 |z_B'(t)|\le C t^{-5/2}\min(t/B,1)
                   +CB^{-1}t^{-3/2}.
\]
Hence $\int_1^\infty|z_B'(t)|\,\mathrm{d}t=O(B^{-1})$, while
$\int_0^1|z_B(t)|\,\mathrm{d}t=O(B^{-1})$. The unit-mesh sum-integral
error is thus $O(B^{-1})$.

In the integral with $a_B=g_B$, substitute $y=g_B(t)$. Then
\[
 t=\frac B\alpha y+\frac{y^2}{2\alpha^2},\qquad
 t^{-3/2}\,\mathrm{d}t=\left(\frac\alpha B\right)^{1/2}
 y^{-3/2}K\!\left(\frac y{\alpha B}\right)\,\mathrm{d}y,
 \quad K(z)=\frac{1+z}{(1+z/2)^{3/2}}.
\]
The linear integral has the same expression with $K$ replaced by 1.
Since $|K(z)-1|\le C\min(z,1)$ and
$|f(q+y)-f(q)|\le C\min(y,1)$, their difference is at most
\[
 CB^{-1/2}\int_0^\infty
 y^{-3/2}\min(y,1)\min(y/B,1)\,\mathrm{d}y=O(B^{-1}).
\]
For the last equality, split at $1$ and $B$. Combining the estimates
proves the assertion.
\end{proof}

After applying Lemma~\ref{curvaturelemma} to~\eqref{renewalapprox}
and using~\eqref{factorremoval}, we may replace $d_m$ by
$\delta-1+\alpha m/\sqrt n$ with an $O(n^{-1/2})$ error.
Canceling $\kappa$ and subtracting
$f(\delta-1)=(\delta-1)^2$ gives the scalar identity
\begin{equation}\label{scalar}
 2\delta-1=\sum_{m\ge1}p_m
 [f(\delta-1+\alpha m/\sqrt n)-f(\delta-1)]
 +O(n^{-1/2}).
\end{equation}

\section{The zeta coefficient}
We now complete the proof of Theorem~\ref{main} by evaluating the
scalar identity~\eqref{scalar}.
Put $\eta=\alpha/\sqrt n$ and $q=\delta-1$. Uniformly for $q$ in a
small fixed neighborhood of $-1/2$,
\[
 |f(q+\eta m)-f(q)|\le C\min(\eta m,1).
\]
By~\eqref{pm}, replacing $p_m$ by
$(2\sqrt\pi)^{-1}m^{-3/2}$ in~\eqref{scalar} costs $O(\eta)$.
Define
\[
 I(q)=\int_0^\infty y^{-3/2}[f(q+y)-f(q)]\,\mathrm{d}y.
\]
A uniform Riemann-sum estimate gives
\begin{equation}\label{Riemann}
 \sum_{m\ge1}m^{-3/2}[f(q+\eta m)-f(q)]
 =\eta^{1/2}I(q)+O(\eta).
\end{equation}
Here are details for the singular endpoint. On a fixed small interval,
the integrand is exactly $2q y^{-1/2}+y^{1/2}$. Its integral over
$(0,\eta)$ is $O(\sqrt\eta)$, and its total variation on
$[\eta,a]$ is $O(\eta^{-1/2})$. On $[a,\infty)$ the profile is
bounded Lipschitz, so the integrand has finite total variation,
uniformly in $q$. The error in the mesh-$\eta$ Riemann sum for the
integral is therefore $O(\sqrt\eta)$; multiplication by
$\eta^{1/2}$ proves~\eqref{Riemann}.

Moreover $I$ is locally Lipschitz near $-1/2$. On a small initial
$y$-interval the difference between its numerators for $q$ and $q_0$
is $2(q-q_0)y$. On the complement use the Lipschitz bound for $f$
and integrability of $y^{-3/2}$. Since
$q+1/2=O(n^{-1/4})$, replacing $I(q)$ by $I(-1/2)$ contributes
$O(n^{-1/2})$ to~\eqref{scalar}.

For $y\ge0$, $f(y-1/2)=h(y-1/2)$, and
\[
 h(y-1/2)-\frac14=\{y\}^2-\{y\}.
\]
Piecewise integration by parts, for an integer $N\ge1$, gives
\begin{align*}
 \int_0^N y^{-3/2}[h(y-1/2)-1/4]\,\mathrm{d}y
 &=8\sum_{k=1}^N\sqrt k-\frac{16}{3}N^{3/2}-4\sqrt N.
\end{align*}
The Euler--Maclaurin expansion
\[
 \sum_{k=1}^N\sqrt k
 =\frac23N^{3/2}+\frac12N^{1/2}+\zeta(-1/2)+O(N^{-1/2})
\]
therefore proves $I(-1/2)=8\zeta(-1/2)$. Substituting
into~\eqref{scalar} yields
\[
 2\delta-1=4\zeta(-1/2)\sqrt{\frac\alpha\pi}\,n^{-1/4}
 +O(n^{-1/2}).
\]
Since $\delta=\alpha\sqrt n-\beta_n$, this is precisely~\eqref{result}
and completes the proof of Theorem~\ref{main}.

\section*{A note on the use of AI}
\label{sec:ai}
\href{https://openai.com/codex/}{OpenAI Codex} (ChatGPT Astra 6) and
\href{https://claude.ai/}{Anthropic Claude} (Fable 5,1) (October 2026) assisted with
numerical code for tables and figures, and final text editing.
The author checked the mathematical arguments and numerical results,
reviewed and edited the text, and takes responsibility for the final
manuscript.

\section*{Funding}
This work was supported by the Czech Science Foundation
grant number 24-11146S.

\section*{Competing interests}
The author declares no competing interests.

\section*{Data availability}
No empirical data were used. The results are established
analytically by the proofs in this paper.

\end{document}